\documentclass[preprint,12pt]{elsarticle}
\usepackage{amssymb, pstricks,pst-plot}
\journal{Comptes Rendus Mathematique}

\newtheorem{theorem}{Theorem}[section]

\newtheorem{corollary}[theorem]{Corollary}

\newtheorem{definition}[theorem]{Definition}

\newcommand{\C}{\mathbb{C}}

\newcommand{\R}{\mathbb{R}}
\newcommand{\T}{\mathrm{T}}

\newcommand{\I}{\mathrm{i}}

\newcommand{\cC}{\mathcal{C}}

\newcommand{\cO}{\mathcal{O}}

\renewcommand{\d}{\mathrm{d}}

\newcommand\wt{\widetilde}

\def\bs{\backslash}

\newcommand\CAP{{\rm CAP}}
\newcommand\PCAP{{\rm PCAP}}

\begin{document}

\begin{frontmatter}
\title{Oka Manifolds}
\author{Franc Forstneri\v c}
\address{Faculty of Mathematics and Physics, University of Ljubljana, 
and Institute of Mathematics, Physics and Mechanics, Jadranska 19, 
1000 Ljubljana, Slovenia}

\begin{abstract}
We give a positive answer to Gromov's 
question in \cite[3.4.(D), p.\ 881]{Gromov:OP}:
If every holomorphic map from a compact convex set in 
a Euclidean space $\C^n$ to a certain complex manifold $Y$ 
is a uniform limit of entire maps $\C^n\to Y$,
then $Y$ enjoys the parametric Oka property. 
In particular, for any reduced Stein space $X$ 
the inclusion $\cO(X,Y) \hookrightarrow \cC(X,Y)$ 
of the space of holomorphic maps into the space of continuous maps
is a weak homotopy equivalence.

\smallskip
\noindent 
{\bf Resume}
\smallskip

\noindent
Nous donnons une r\'eponse positive \`a la question suivante pos\'ee par 
Gromov dans \cite[3.4.(D), p.\ 881]{Gromov:OP}:
Si une vari\'et\'e analytique complexe $Y$ est telle que toute application 
holomorphe d'un voisinage d'un sous-ensemble compact convexe de l'espace 
euclidien $\C^n$ dans $Y$ peut \^etre approxim\'ee par des applications 
enti\`ere de $\C^n$ dans $Y$, alors les applications holomorphes d'un espace 
de Stein r\'eduit $X$ dans $Y$ poss\`edent la propri\'et\'e d'Oka 
param\'etrique. 
\end{abstract}


\end{frontmatter}

\medskip
\centerline{\em To Mikhael L.\ Gromov on the occasion of receiving the Abel Prize}

\section{The Oka-Grauert-Gromov Principle}
\label{Sec:Oka}
A complex manifold $Y$ is said to enjoy the 
{\em Convex Approximation Property} (CAP)
if every holomorphic map from a neighborhood of a compact convex 
set $K\subset \C^n$ to $Y$ can be approximated, uniformly on $K$, 
by entire maps $\C^n\to Y$. This property was introduced in
\cite{FF:CAP,FF:EOP} where it was shown to 
be equivalent to the basic Oka property of $Y$.

Here we show that CAP also implies the {\em parametric Oka property},
thereby providing a positive answer to Gromov's question 
\cite[3.4.(D), p.\ 881]{Gromov:OP}.  Our main result is the following:

\begin{theorem}
\label{CAP-POP}
If $h\colon Z\to X$ is a stratified holomorphic fiber bundle 
over a reduced Stein space $X$ such that all its fibers satisfy  \CAP,
then sections $X\to Z$ satisfy the parametric Oka principle. 
In particular, the inclusion of the space of 
all holomorphic sections $X\to Z$ into the space
of all continuous sections is a weak homotopy equivalence.
\end{theorem}

The conclusion of Theorem \ref{CAP-POP} means the following:

\smallskip
\noindent 
{\em Given a compact $\cO(X)$-convex subset $K$ of $X$, 
a closed complex subvariety $X'$ of $X$, compact sets
$P_0\subset P$ in a Euclidean space $\R^m$,
and a continuous map $f\colon P\times X\to Z$ such that 
(a)  for every $p\in P$, $f(p,\cdotp)\colon X\to Z$ is a section 
of $Z\to X$ that is holomorphic on a neighborhood of $K$
(independent of $p$) and such that $f(p,\cdotp)|_{X'}$ is holomorphic on $X'$, and 
(b) $f(p,\cdotp)$ is holomorphic on $X$ for every $p\in P_0$,
then there is a homotopy $f_t\colon P\times X \to Z$ $(t\in [0,1])$, 
with $f_0=f$, such that $f_t$ enjoys properties (a) and (b) 
for all $t\in[0,1]$, and
\begin{itemize} 
\item[(i)]    $f_1(p,\cdotp)$ is holomorphic on $X$ for all $p\in P$, 
\item[(ii)]   $f_t$ is uniformly close to $f$ on $P\times K$ for all $t\in [0,1]$, and 
\item[(iii)]  $f_t=f$ on $(P_0\times X)\cup (P\times X')$ for all $t\in [0,1]$.
\end{itemize}
If in addition $f(p,\cdotp)$ is holomorphic on a neighborhood of 
$X'$ for all $p\in P$, then the homotopy $f_t$ can be chosen fixed 
to a given finite order along $X'$.
}
\smallskip

In the special case when $Z\to X$ is a holomorphic fiber bundle whose fiber 
is a complex homogeneous manifold, Theorem \ref{CAP-POP} 
is due to Grauert \cite{Grauert3}. 
In 1989 Mikhael Gromov published an influential paper \cite{Gromov:OP}
in which he obtained Theorem \ref{CAP-POP} when $Z\to X$ is a fiber bundle over 
a Stein manifold $X$ whose fiber $Y$ is {\em elliptic} in the sense that it admits a 
{\em dominating holomorphic spray} --- a triple $(E,\pi,s)$ consisting 
of a holomorphic vector bundle $\pi \colon E\to Y$ 
and a holomorphic map $s\colon E\to Y$ 
such that for each $y\in Y$ we have $s(0_y)=y$
and $(\d s)_{0_y}(E_y)=\T_y Y$. 
It is easily seen that ellipticity
implies CAP, but the converse implication is not known.  
In the same paper Gromov asked whether a Runge type approximation property for 
holomorphic maps $\C^n \to Y$ might suffice to infer the Oka principle
\cite[3.4.(D), p.\ 881]{Gromov:OP}. 
Theorem \ref{CAP-POP} gives a positive answer to Gromov's question
and shows that all Oka properties considered in the literature are equivalent 
to each other. 
Hence this is an opportune moment to introduce 
the following class of complex manifolds.

\begin{definition}
A complex manifold $Y$ is said to be an {\em Oka manifold} if
it enjoys \CAP\ or, equivalently, the parametric Oka property 
for maps of all reduced Stein spaces to $Y$.
\end{definition}

Theorem \ref{CAP-POP} indicates that the class of Oka manifolds 
is dual to the class of Stein manifold in a sense that can be made 
precise by means of abstract homotopy theory 
(see L\'arusson \cite{Larusson2,Larusson3}).

Examples of Oka manifolds can be found in the papers 
\cite{Gromov:OP,FF:CAP,FF:EOP}.
The following result is useful for finding new examples.

\begin{corollary}
\label{ascend-descend}
{\rm (C.f.\ \cite[Theorem 1.4]{FF:CAP})}
Assume that $E$ and $B$ are complex manifolds. If $\pi \colon E\to B$ 
is a holomorphic fiber bundle whose fiber $Y$ is an Oka manifold,
then $E$ is an Oka manifold if and only if $B$ is an Oka manifold.
\end{corollary}

\section{A characterization of the parametric Oka property}
\label{known}
We recall from \cite{FF:CAP} the precise definition of CAP
and its parametric analogue, PCAP. 
Let $z=(z_1,\ldots,z_n)$ be complex coordinates on $\C^n$, with
$z_j=x_j+\I y_j$. A {\em special convex pair} $(K,L)$ in $\C^n$ consists 
of a cube 
$    
		L =\{z\in \C^n\colon |x_j| \le a_j,\ |y_j|\le b_j,\ j=1,\ldots,n\}
$
and a compact convex set 
$
    K =\{z\in L \colon y_n \le h(z_1,\ldots,z_{n-1},x_n)\},
$
where $h$ is a continuous concave function with values
in $(-b_n,b_n)$. 

We say that a map is holomorphic on a compact set if it is holomorphic
in an open neighborhood of that set.

\begin{definition}
\label{def:CAP}
A complex manifold $Y$ enjoys \CAP\ if for each special convex pair $(K,L)$ in $\C^n$, 
every holomorphic map $f\colon K\to Y$ can be approximated uniformly on $K$ 
by holomorphic maps $L \to Y$.

$Y$ enjoys the {\em Parametric Convex Approximation Property} {\rm (PCAP)}
for a pair of topological spaces $P_0\subset P$ if for every special 
convex pair $(K,L)$, a continuous map $f\colon P\times L \to Y$ 
such that $f(p,\cdotp)\colon L\to Y$ is holomorphic for every $p\in P_0$,
and is holomorphic on $K$ for every $p\in P$, can be approximated 
uniformly on $P\times K$ by continuous maps $\wt f\colon P\times L \to Y$
such that $\wt f(p,\cdotp)$ is holomorphic on $L$ for 
all $p\in P$, and $\wt f=f$ on $P_0\times L$.
\end{definition}

The following result was proved in \cite{FF:CAP,FF:EOP,FF:Kohn}.

\begin{theorem} 
\label{CAP}
(a) If $Y$ enjoys \CAP, then it enjoys the basic Oka property 
(the conclusion of Theorem \ref{CAP-POP} for $P$ a singleton and $P_0=\emptyset$).

(b) If \ $Y$ enjoys \PCAP\ for a  pair of compact Hausdorff 
spaces $P_0\subset P$, then it also satisfies the conclusion of
Theorem \ref{CAP-POP} for the same pair.
\end{theorem}

Hence Theorem \ref{CAP-POP} follows from the following result.

\begin{theorem}
\label{CAP-PCAP}
If $Y$ enjoys {\rm CAP}, then it also enjoys {\rm PCAP} for all pairs 
$P_0\subset P$ of compact parameter spaces contained in a 
Euclidean space $\R^m$.
\end{theorem}

\section{Proof of Theorem \ref{CAP-PCAP}}
Assume that $P_0\subset P$ are compact sets in $\R^m \subset \C^m$,
$(K,L)$ is a special convex pair in $\C^n$, 
$U\supset K$ and $V\supset L$ are open convex neighborhoods 
of $K$ resp.\ $L$ in $\C^n$, and $f\colon P\times V \to Y$ 
is such that $f(p,\cdotp)\colon V \to Y$ 
is holomorphic for every $p\in P_0$, and it is holomorphic on 
$U$ for every $p\in P$. 

By using a Tietze extension theorem for Banach-valued maps
as in \cite[Proposition 4.4]{FF:Rothschild} and shrinking 
the set $V \supset L$ we can assume that $f(p,\cdotp)$ is holomorphic on 
$V$ for all $p$ in a neighborhood $P'_0\subset \C^m$ of $P_0$.

We may assume that $0\in \C^n$ belongs to ${\rm Int} K$.
Choose a continuous function $\tau\colon P\to [0,1]$ 
such that $\tau=0$ on $P_0$ and $\tau=1$ on $P\bs P'_0$. Set
\[
		f_t(p,z) = f\bigl(p,(1-(1-t)\tau(p))z\bigr)
		\in Y, \quad p\in P,\ z\in V,\ t\in[0,1].
\]	 
Then $f_t$ has the same properties as $f=f_1$, the homotopy is fixed 
for $p\in P_0$, and the map $f_0(p,\cdotp)$ is holomorphic on
$V$ for all $p\in P$.

We shall follow Gromov's  proof of the h-Runge theorem in \cite{Gromov:OP},
using local sprays over Stein domains covering the graph
of the homotopy $f_t$, along with the basic Oka property of $Y$
(which is implied by CAP in view of Theorem \ref{CAP} (a)). Set 
$Z=\C^m\times \C^n\times Y$. For every  $t\in [0,1]$ let
\[
		F_t(p,z)= \left( p,z,f_t(p,z)\right), \quad
		\Sigma_t = F_t(P\times K) \subset Z.
\]
We also set $S_0= F_0(P\times L)$. 
By \cite[Corollary 2.2]{F-Wold} the sets 
$\Sigma_t$ and $S_0$ are Stein compacta in $Z$. 
Hence there are numbers $0=t_0<t_1<\ldots <t_N=1$
and Stein domains $\Omega_j\subset Z$ such that 
\begin{equation}
\label{Omegaj}
		\Sigma_t \subset \Omega_j \ {\rm when}\  
		t_{j}\le t\le t_{j+1}\ {\rm and}\ j=0,1,\ldots,N-1.
\end{equation}

Let $\pi_Y \colon Z  \to Y$ denote the projection onto $Y$, 
and let $E=\pi_Y^* (TY)$ denote the pull-back of the tangent bundle 
of $Y$ to $Z$. By standard techniques we obtain for every Stein domain 
$\Omega\subset Z$ a Stein neighborhood $W \subset E|_\Omega$ of the zero 
section $\Omega \subset E|_\Omega$ and a holomorphic map $s\colon W\to Z$ 
that maps the fiber $W_{(\zeta,z,y)}$ over $(\zeta,z,y)\in Z$
biholomorphically onto a neighborhood of this point in 
$\{(\zeta,z)\} \times Y$ such that $s$ preserves the
zero section $\Omega$ of $E|_\Omega$. In short, $s$ 
is a local fiber dominating spray in the sense of 
Gromov \cite{Gromov:OP}. We may assume that $W$ is 
Runge in $E|_\Omega$ and its fibers are convex (in the fibers of $E$). 
Since  $E|_{\Omega}$ is Stein and 
$Y$ enjoys \CAP\ (and hence the basic Oka property by 
Theorem \ref{CAP} (a)), $s$ can be approximated uniformly
on compacts in $W$ by a global spray $\wt s\colon E|_{\Omega}\to Z$ 
that agrees with $s$ to the second order along 
the zero section $\Omega$. 

\begin{figure}[ht]
\label{F1}
\psset{unit=0.5cm,linewidth=0.7pt}

\begin{pspicture}(-8.3,-0.5)(12,10)

\pspolygon(-3,0)(8,0)(12,6)(1,6)             
\psline[linewidth=0.4pt](-1,3)(10,3)         

\psline(-5,0)(-1,6)
\rput(-4.4,3){$(p,z)$}
\rput(-6,1.2){$\C^m\times \C^n$}

\psdot[dotsize=3pt](-3,3)
\psline[linewidth=0.2pt]{<-}(-2.8,3)(-1.2,3)

\pscurve[linewidth=1pt](0.5,0)(1,0.5)(3,2.5)(5,5)(5.8,6)       
\psecurve[linewidth=1pt](0.5,0)(1,0.5)(3,2.5)(5,5)(5.8,6)      
\psecurve[linewidth=1pt](2.5,0)(3,0.5)(5,2.5)(6.9,5)(7.5,6)    

\pscustom[fillstyle=solid,fillcolor=lightgray]      
{
\pscurve(0.5,0)(1,0.5)(3,2.5)(5,5)(5.8,6)
\psline(5.8,6)(5.8,10)
\pscurve(5.8,10)(5,9.3)(3,7)(1,5)(0.5,4.5)
\psline(0.5,4.5)(0.5,0)
}

\psecurve[linewidth=1pt](0.5,1.8)(1,2.3)(3,4)(5,6.7)(5.8,7.5)  
\rput(2.7,5){$\xi_t$}
\psline[linewidth=0.2pt]{->}(2.7,4.7)(2.7,3.8)

\rput(10.5,4.7){$Z_{(p,z)} = Y$}
\psline[linewidth=0.2pt]{->}(9.5,4.2)(9.5,3.05)

\rput(5.5,1.2){$F_t$}
\psline[linewidth=0.2pt]{<-}(3.8,1.2)(4.8,1.2)

\rput(9,7.8){$F_0(P\times V)$}
\psline[linewidth=0.2pt]{<-}(5.5,5.5)(8.2,7.3)

\rput(1,8.1){$E|_{F_0(P\times V)}$}
\psline[linewidth=0.2pt]{->}(1.3,7.5)(2.5,6.5)

\psarc[linewidth=0.2pt,arrows=<-](2,2){4}{21}{58}     
\rput(5.5,4.5){$s$}

\rput(8,1.2){$Z$}

\end{pspicture}
\caption{Lifting sections $F_t$ to the spray bundle $E|_{F_0(P\times V)}$}
\label{Fig:1}
\end{figure}
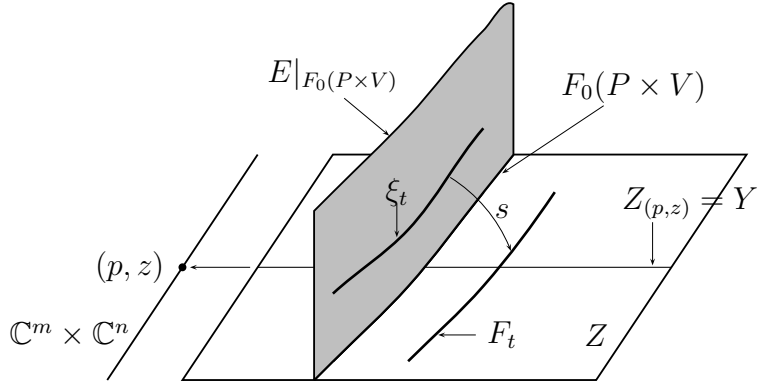

This shows that after refining our subdivision 
$\{t_j\}$ of $[0,1]$ and shrinking the set $U\supset K$
there are Stein domains $\Omega_0,\ldots,\Omega_{N-1}$ 
as in (\ref{Omegaj}), sprays $s_j \colon E|_{\Omega_j}\to Z$, 
and homotopies of $z$-holomorphic sections $\xi_t$ 
$(t\in [t_{j},t_{j+1}])$ of the restricted bundle 
$E|_{F_{t_j}(P\times U)}$ 
such that $\xi_{t_j}$ is the zero section,  
$\xi_{t}(p,\cdotp)$ is independent of $t$ when $p\in P_0$ 
(hence it is the zero section), and  
\[ 
		s_j \circ \xi_t \circ F_{t_{j}} = F_t 
		\quad {\rm on} \ P\times U, \quad t\in [t_j,t_{j+1}].
\]
Furthermore, the existence of such liftings $\xi_t$ 
is stable under sufficiently small perturbations of 
the homotopy $F_{t}$. (See Fig.\ \ref{Fig:1}.)

%
%
%
%

Consider the homotopy of sections $\xi_t$ of 
$E|_{F_0(P\times U)}$ for $t\in [0,t_1]$.
By the Oka-Weil theorem we can approximate 
$\xi_t$ uniformly on $P\times K$ by $z$-holomorphic sections 
$\wt \xi_t$ of $E|_{F_0(P\times V')}$ for an open convex set
$V'\subset \C^n$ with $L\subset V'\subset V$. (This parametric
version of the Oka-Weil theorem is obtained by using 
a continuous partition of unity with respect to the parameter.)
Further, we may choose $\wt\xi_t=\xi_t$ 
for $t=0$ and on $P_0\times V'$.

By \cite[Corollary 2.2]{F-Wold} there is a  Stein neighborhood 
$\Omega \subset Z$ of $S_0$ such that $\Sigma_0$ is 
$\cO(\Omega)$-convex. Hence $E|_{\Sigma_0}$  is exhausted 
by $\cO(E|_\Omega)$-convex compact sets. Since $E|_{\Omega}$ is Stein and 
$Y$ enjoys \CAP, the spray $s$ can be approximated on the range of 
the homotopy $\{\xi_t\colon t\in [0,t_1]\}$ 
by a spray $\wt s\colon E|_{\Omega}\to Z$ that 
agrees with $s$ to the second order along the zero section. 
The maps
\[
		\wt f_t = \pi_Y\circ \wt s \circ \wt \xi_t\circ F_{0}
		\colon P\times V' \to Y ,\quad t\in [0,t_1]
\]
are then $z$-holomorphic on $V'\supset L$, and they
approximate $f_t$ uniformly on $P\times K$.
If the approximation is sufficiently close, 
we obtain a new homotopy $\{f_t\colon t\in [0,1]\}$ that agrees 
with $\wt f_t$ for $t\in [0,t_1]$ (hence is $z$-holomorphic on
$L$ for these values of $t$), and that agrees with the initial
homotopy for $t\in  [t'_1,1]$ for some $t'_1>t_1$ close to $t_1$.

We now repeat the same argument with the parameter interval 
$[t_1,t_2]$, using $f_{t_1}$ as the new reference map 
(that is $z$-holomorphic over $L$). This gives us a new homotopy 
that is $z$-holomorphic on $L$ for $t\in [0,t_2]$. 

After finitely many steps of this kind we obtain a desired homotopy 
whose final map $\wt f$ at $t=1$ satisfies the conclusion of 
Theorem \ref{CAP-PCAP}.


\begin{thebibliography}{00}

\bibitem{F-Wold}
Forstneri\v c, F., Wold, E.\ F.:
Fibrations and Stein Neighborhoods.
Preprint (2009). arXiv: 0906.2424 


\bibitem{FF:CAP}
Forstneri\v c, F.:
Runge approximation on convex sets implies Oka's property.
Ann.\ Math.\ (2) \textbf{163}, 689--707 (2006)

\bibitem{FF:EOP}
Forstneri\v c, F.:
Extending holomorphic mappings from subvarieties in Stein manifolds.
Ann.\ Inst.\ Fourier  \textbf{55}, 733--751 (2005)

\bibitem{FF:Kohn}
Forstneri\v c, F.:
The Oka principle for sections of stratified fiber bundles.
Pure and Appl.\ Math.\ Quarterly, (2009).
arXiv: 0705.0591


\bibitem{FF:Rothschild}
Forstneri\v c, F.:
Invariance of the parametric Oka property.
(Proceedings of the conference in honor of Linda P.\ Rothschild,
Fribourg, Switzerland, July 2008), Birkh\"auser Verlag, to appear.
arXiv: 0705.0591

\bibitem{Grauert3}
Grauert, H.:
Analytische Faserungen \"uber holomorph-vollst\"andigen R\"aumen.
Ma\-th.\ Ann.\ \textbf{135},  263--273 (1958)

\bibitem{Gromov:OP}
Gromov, M.:
Oka's principle for holomorphic sections of elliptic bundles.
J.\ Amer.\ Math.\ Soc.\ \textbf{2}, 851-897 (1989)

\bibitem{Larusson2}
L\'arusson, F.:
Model structures and the Oka principle.
J.\ Pure Appl.\ Algebra\ \textbf{192}, 203--223 (2004)

\bibitem{Larusson3}
L\'arusson, F.:
Mapping cylinders and the Oka principle.
Indiana Univ.\ Math.\ J.\ \textbf{54}, 1145--1159 (2005)


\end{thebibliography}
\end{document}